\documentclass[12pt]{article}
\usepackage{graphicx}
\usepackage{amssymb}
\usepackage{amsmath}
\usepackage{amsthm}
\usepackage{cite}
\newcommand{\norm}[1]{\left\Vert#1\right\Vert}

\numberwithin{equation}{section}
\def\nn{{\nonumber}}

\def\BR{{\mathbb R}}
\def\BN{{\mathbb N}}
\def\BZ{{\mathbb Z}}

\def\BT{{\mathbb T}}
\newtheorem{Pa}{Paper}[section]
\newtheorem{Tm}[Pa]{{\bf Theorem}}
\newtheorem{La}[Pa]{{\bf Lemma}}
\newtheorem{Cy}[Pa]{{\bf Corollary}}
\newtheorem{Rk}[Pa]{{\bf Remark}}

\newtheorem{Ee}[Pa]{{\bf Example}}
\newtheorem{Dn}[Pa]{{\bf Definition}}
\newtheorem{Pn}[Pa]{{\bf Proposition}}

\def\XXint#1#2#3{{\setbox0=\hbox{$#1{#2#3}{\int}$}
     \vcenter{\hbox{$#2#3$}}\kern-.5\wd0}}

\title{Differential equations on a $k$-dimensional torus: Poincar{\'e}  type results}
\author{Lev Sakhnovich}
\date{}
\begin{document}

\maketitle

\begin{abstract}
Ordinary differential equations of the first order on the torus have been investigated in detail by H. Poincar{\'e} and A. Denjoy.
The long-standing problem of generalising these results for the equations of the order $k>1$ (or for the systems of equations) 
is both important and  difficult, requiring an essentially new approach.
(See, e.g., an interesting old paper by P.~Bohl.) 

In this paper, we propose a new (non-Hamiltonian) and promising approach.
 We  use Hamiltonians, that is,   ordinary differential systems of equations of the first order, only for heuristics.
In the main   scheme  and corresponding proofs we do not use these systems.\
Instead  of differential systems, we study   sets of continuous vector functions $\phi(t,\eta)$ satisfying   important conditions,
which follow from the analogy with the solutions in the case $k=1$.
Some of  our results are new even in the case $k=1.$

\end{abstract}


\vspace{0.3em}


\section{Introduction}
Ordinary differential equations of the first order on the torus have been investigated in detail by H. Poincar{\' e}  and A. Denjoy.
P. Bohl,  back in 1916 (see \cite{Bo}), emphasised the importance  of the transfer of the results for the order $k=1$ to the case $k>1$, adding at the same time:
``However, any attempt to do so would be hopeless". It became clear that a new approach to this problem is needed.

In this paper, we propose a new (non-Hamiltonian) and promising approach.
 We  use Hamiltonians, that is,   ordinary differential systems of equations of the first order, only for heuristics.
In the main   scheme  and corresponding proofs we do not use these systems.
Instead  of differential systems, we study   sets of continuous vector functions $\phi(t,\eta)$ satisfying   important conditions,
which follow from the analogy with the solutions in the case $k=1$.
We note that a shift from  differential equations is not new in physics
(see, e.g., \cite{Fey} and our Remark \ref{Remark 1.4}).

In a different context, rotation vectors and rotation sets (limit sets) appear in the interesting works \cite{ARS, Fra1, Lel, Lli0, Lli}
(see also the references therein).

Let us formulate some classical Poincar{\'e} results \cite[Ch. XVII]{CL}. For this purpose,
consider the ordinary differential equation
\begin{equation}\frac{d}{dt}{x(t)}=f(t,x),\label{1.1}\end{equation}
where the function $f(t,x)$  is  real, continuous and defined for all real numbers $t$ and real numbers  $x$.
It is assumed in this theory that two conditions below are satisfied.
\begin{equation}{\mathbf{Condition \, I.}}\qquad \qquad \qquad \qquad f(t+1,x)=f(t,x+1)=f(t,x). \label{1.2} \end{equation}
\textbf{\quad Condition II.} \emph{Through every point of the $(t,s)$ plane passes a unique solution of equation \eqref{1.1}.}

Let us  study  solutions of the differential system \eqref{1.1} with initial value (parameter) $\eta$.
It is assumed  that the function $\phi(t,\eta)$ is known at the two fixed time points
$t=0 $ and $t=1$:
\begin{equation}x=\phi(t,\eta), \quad \phi(0,\eta)=\eta, \quad \phi(1,\eta)=\psi(\eta).\label{1.3}\end{equation}
It follows from \eqref{1.1}--\eqref{1.3}
that
\begin{equation}\phi(t,\eta+1)=\phi(t,\eta)+1,\quad \psi(\eta+1)=\psi(\eta)+1. \label{1.4}\end{equation}
One may easily see that
\begin{equation}\phi(t+1,\eta)=\phi(t,\psi(\eta)). \label{1.5}\end{equation}
Indeed, in view of \eqref{1.2} the function $\phi(t+1,\eta)$ is a solution of \eqref{1.1} because $\phi(t,\eta)$ is its solution.
Moreover, it follows from \eqref{1.3} that $\phi(1,\eta)=\phi(0,\psi(\eta))=\psi(\eta)$ and so the solutions
$\phi(t+1,\eta)$ and $\phi(t,\psi(\eta))$ both have the same both have the same value $\psi(\eta)$ at $t=0$.
Hence,  \eqref{1.5} is valid for all $t$ (compare with \cite[Ch.1, Sect. 1]{CL}).
We use the uniqueness theorem, and not the condition II, since the uniqueness theorem is also true in the vector case (see \cite[Ch.~II, Sect. 1]{Lef}), which we consider 
in the next section.
\begin{Rk}\label{Remark 1.1}Two points  $P_1=(t_1,x_1)$ and  $P_2=(t_2,x_2)$  are regarded
as identical if  $t_1-t_2$ and $x_1-x_2$  are integers.\end{Rk}
One may represent the solution paths $(t,x)$ on torus in the following way:
\begin{align}& u=\big(a+b\cos(2\pi{x})\big)\cos(2\pi{t}),\label{1.6}
\\ &
v=\big(a+b\cos(2\pi{x})\big)\sin(2\pi{t}),\label{1.7}
\\ &
w=b\sin(2\pi{x}),\label{1.8}
\end{align}
where $a$ and $b$ are constants such that $0<b<a$, and $(u,v,w)$ are rectilinear coordinates  in a three-dimensional space.
Now, let $x=\phi(t,\eta)$ be the solution of \eqref{1.1} such that $\phi(0,\eta)=\eta$ (see \eqref{1.3}).
Let us consider the function
\begin{equation}\psi(\eta)=\phi(1,\eta)\label{1.9} .\end{equation}
According to  \cite[Ch. XVII, Sect. 1]{CL}, the function $\psi(\eta)$ is a continuous, monotonous  and increasing homeomorphism of the real line into itself and (see \eqref{1.4})
\begin{equation}\psi(\eta+1)=\psi(\eta)+1.\label{1.10}\end{equation}
The function $\psi(\eta)$ \emph{represents} a transformation $T$  of the form $TP=P_1$  or, introducing it in a slightly different way,
\begin{equation}T\eta=\psi(\eta).\label{1.11}\end{equation}
\begin {Tm}\label{Theorem 1.2} \cite[Ch. XVII, Sect. 2, Theorem 1.2]{CL} The limit
\begin{equation}\rho=\lim_{t{\to}\infty}\frac{\phi(t,\eta)}{t}\label{1.12} \end{equation}
exists and does not depend on the initial value $\eta$; it is rational if and only if some integer and positive power $m\in \BN$ of $T$ has a fixed point.\end{Tm}
Recall that $\BN$ stands for the set of positive integers.
\begin{Dn}\label{Definition 1.3} The value $\rho$ introduced in \eqref{1.12} is called  the rotation number.
\end{Dn}
Let $\psi_{n}(\eta)$ be the function defined by the relations
\begin{equation}
\psi_{0}(\eta)=\eta, \quad \psi_{n}(\eta)=\psi[\psi_{n-1}(\eta)] \quad (n\in\BN).
\label{1.14}\end{equation}
 The  function $\psi_{n}$ is of the same type as $\psi$, that is,   $\psi_n(\eta)$ is continuous, monotonous and satisfies \eqref{1.10}.
Next, $\rho$ is expressed in terms of $\psi_n$.
\begin {Tm}\label{Theorem 1.4}  \cite[Ch. XVII, Sect. 5, Theorem 1.3]{CL} The rotation number satisfies the equality
\begin{equation}\rho=\lim_{n{\to}\infty}{\psi_n}(\eta)\big/{n},\label{1.14+} \end{equation}
where the limit exists and does not depend on the initial value $\eta$. \end{Tm}

In this paper, we propose a new (non-Hamiltonian) and promising approach.  Further we consider the $k$-dimensional vector space
$E_k$ of vectors $X=[x_1,x_2,...,x_k]$ with real-valued elements, where the norm   is defined by the relation
\begin{equation}\norm{X}=\max|x_p|,\quad 1{\leq}p{\leq}k.\label{1.15}\end{equation}
\begin{Dn}\label{Definition 1.5}$($see\cite{Vul}$)$ Let vectors $X$ and $Y$belong to $E_k$. We say that $X>Y$ if all the elements of the
vector $X-Y$ are non-negative and at least one of them is positive.\end{Dn}
Thus, $E_k$ is a semi-ordered space.
\vspace{0.5em}

{\it Our approach has the following characteristic features:}

1. We  use  Hamiltonians, that is,   ordinary differential systems of equations of the first order of the form  \eqref{1.1} only as a source of heuristics.
The corresponding differential systems are  not used in our statements and proofs.

2. Instead of  differential systems, we study a set of continuous vector functions $\phi(t,\eta)$, which satisfy  vector versions of conditions  \eqref{1.3}--\eqref{1.5}
and which we call solutions of {\it generalised  vector systems} (see section 2).
We do not assume that the functions $\phi(t,\eta)$ are differentiable.

3. We consider the corresponding problems in  $k$-spaces $E_k$  $(k{\geq}1)$.

4. We use the theory of the semi-ordered spaces  $E_k$ and Zorn's lemma.

5. Under natural conditions, we construct solutions $\phi(t,\eta)$ in an explicit form (see Theorem \ref{Theorem 2.3}).
This result is new even for the case $k=1$.

6. We obtain   vector versions  of Theorems  \ref{Theorem 1.2} and   \ref{Theorem 1.4} (see section 3).

7. It is very difficult  to check the critical points of the solution $\phi(t,\eta)$ within the framework of the classical theory.  Our approach allows to  effectively solve
this problem (see Section 4).
 \begin{Rk}\label{Remark 1.4}  The avoidance  of differential equations is not new in physics.  In particular, the overcomplexity of the
quantum electrodynamics  equations is well known. R. Feynman \cite{Fey}  built an effective and fairly simple procedure for solving a number of problems of
quantum electrodynamics without using the corresponding  differential equations. Sure, classical equations served as a background and a hint for heuristics.
Interesting interconnections between classical theory and R. Feynman theory are studied in our paper \cite{Sak}.\end{Rk}

The notation $\BZ$ stands for integer and the notation $\BR$  for  real numbers.
\section{Generalised  vector systems on torus}
\subsection{Differential systems}
We start with a system of differential equations on a torus:
\begin{equation}\frac{d}{dt}{x(t)}=f(t,x),\label{2.1}\end{equation}
where the elements of the $k$-dimensional vector  function $f(t,x)$  are  real, continuous and defined for all numbers $t \geq 0$ and real vectors  $x(t)\in E_k$.

{\it Further we always assume that $t \geq 0$. The case $t \leq 0$ is easily   reduced  to the case $t \geq 0$.}

It is assumed that
\begin{equation} f(t+1,x)=f(t,x+q)=f(t,x), \label{2.2} \end{equation}
for each  $k$-dimensional vector $q$ with integer elements (i.e., for $q\in \BZ^k$).

The vector solution of the differential system \eqref{2.1} is denoted by $\break x=\phi(t,\eta),$
and we assume again that the  function  $\phi(t,\eta)$ is known at two fixed
time points $t=0$ and $t=1$:
\begin{equation}\phi(0,\eta)=\eta, \quad \phi(1,\eta)=\psi(\eta)  \quad (\eta \in E_k).\label{2.3}\end{equation}

\begin{Pn}\label{Proposition 2.1} Let $x=\phi(t,\eta)$ satisfy \eqref{2.1} and \eqref{2.2}.

Then, we have
\begin{align}& \phi(t,\eta+q)=\phi(t,\eta)+q,\quad \psi(\eta+q)=\psi(\eta)+q \quad {\mathrm{for}} \quad q\in \BZ^k; \label{2.4}\\
&
\phi(t+1,\eta)=\phi(t,\psi(\eta)).\label{2.5}\end{align}\end{Pn}
\begin{Rk}\label{Remark 2.2} We note that relations \eqref{2.4} and \eqref{2.5} are vector versions of \eqref{1.4} and \eqref{1.5}
and are proved in the same way as \eqref{1.4} and \eqref{1.5}.
\end{Rk}
It follows from \eqref{2.4} and  \eqref{2.5} that
\begin{align}& \phi(t,\eta)=\eta+a(t,\eta) \label{2.6}
\\ &
\psi(\eta)=\eta+a(\eta),\label{2.7} \end{align}
where
\begin{equation}a(t,\eta)=a(t,\eta+q) \quad {\mathrm{for}} \quad q\in \BZ^k,\quad a(\eta)=a(1,\eta) .\label{2.8}\end{equation}

\subsection{Generalised systems}
Neither  equation \eqref{2.1} nor condition \eqref {2.2} are used in our further considerations.
{\it We   consider a set of real, continuous
$k$-dimensional vector  functions $\phi(t,\eta)$, where $t\geq 0$ and  $\eta\in E_k$.
We assume that   relations \eqref{2.3}--\eqref{2.5} are fulfilled.}
Such a set we call a {\it generalised
system}.

Here and further,  the continuity of the functions depending $t$ and $\eta$
means that they continuously  act on the elements $(t,\eta)$ of the set $\{(t,\eta)\}$.
The requirement for the function $\phi(t,\eta)$ to be differentiable is omitted.
The restriction of $\phi$ on the interval $[0,1]$ is denoted by $\Phi$:
\begin{equation}\Phi(t,\eta)=\phi(t,\eta)\quad (0 \leq t \leq 1,\quad \eta \in E_k). \label{2.9}\end{equation}
\begin{Tm}\label{Theorem 2.3}
{\rm I.} Let a continuous vector function $\phi(t,\eta)$ satisfy
 relations \eqref{2.3}-\eqref{2.5}.
Then, the continuous vector function $\Phi(t,\eta)$
satisfies the following conditions.

a) The relation
\begin{equation}\Phi(t,\eta+q)=\Phi(t,\eta)+q,\label{2.10}\end{equation}
holds for the   vectors $q$ with integer elements $($i.e., for $q\in \BZ^k)$.

b) The following equalities are valid:
\begin{equation}\Phi(0,\eta)=\eta , \quad \Phi(1,\eta)=\psi(\eta).\label{2.11}\end{equation}

{\rm II.} Let the continuous vector function
$\Phi(t,\eta)$  $(0\leq t \leq1,\,\, \eta \in E_k)$ satisfy  conditions  a) and b).
Then,  the solution $\phi(t,\eta)$ of the generalised system is recovered  in the explicit form
via the formula
\begin{equation}\phi(t+n,\eta)=\Phi(t,\psi_n(\eta)) \quad (0 \leq t \leq 1), \label{2.12}\end{equation}
where $n\in \BN$ and $\psi_n$ is expressed via $\psi$ in \eqref{1.4}.
This  $\phi(t,\eta)$ is continuous and satisfies    relations \eqref{2.3}--\eqref{2.5}.
\end{Tm}
\begin{proof}  Part I of the theorem immediately follows from the theorem's conditions.
Note that  formula \eqref{2.12} follows  (under Part I conditions)  from \eqref{2.5} and \eqref{2.9}.

Let us prove Part II.
First, given \eqref{2.12}, let us prove that the constructed vector function $\phi(t,\eta)$ is continuous.
Indeed, discontinuities may only occur at the points $t=n$ for $n\in \BN$.
(Recall that the case $t\geq 0$ is considered.)
In view of \eqref{2.11} and \eqref{2.12}, we have the equalities:
\begin{align}& \phi(n+0,\eta)=\Phi(+0,\psi_n(\eta))=\Phi(0,\psi_n(\eta))=\psi_n(\eta), \label{2.13}
\\ &
\phi(n-0,\eta)=\Phi(1-0,\psi_{n-1}(\eta))=\Phi(1,\psi_{n-1}(\eta))=\psi_n(\eta). \label{2.14}\end{align}
Thus, the constructed vector function $\phi(t,\eta)$ is continuous. Now, formula \eqref{2.5} follows from \eqref{2.9} and \eqref{2.12}.
Relations \eqref{2.10} and \eqref{2.11} yield the second equality in \eqref{2.4}. The second equality in \eqref{2.4} and relations \eqref{2.5}, \eqref{2.9} and \eqref{2.10}
imply the first equality in \eqref{2.4}. Formulas \eqref{2.3} follow from \eqref{2.9} and \eqref{2.10}.  In this way, we proved that the constructed solution $\phi(t,\eta)$ satisfies \eqref{2.3}--\eqref{2.5}.\end{proof}

Theorem \ref{Theorem 2.3} is new even in the scalar case.

\section{ Rotation vectors}
Recall the definition of the space $E_k$ with the norm \eqref{1.15} and a partial order on $E_k$ given by Definition \ref{Definition 1.5} in introduction. 
 \begin{Rk}\label{Remark 3.1} Similar to the case $k=1$  in Remark \ref{Remark 1.1}, the points  $\break P_1=(t_1,X)$ and  $P_2=(t_2,Y)$  $(X,Y \in E_k)$ are regarded
as identical if  $t_1-t_2$ and all the elements of the vector $X-Y$  are integers. This condition shows that differential equations on the torus ${\bf T}^{k}$
are equivalent to corresponding differential equations in the space $E_k$.\end{Rk}
 
Note that formula \eqref{2.7} coincides with Poincare mapping in the case $k=1$.
Similar to \eqref{1.11}, we introduce a transformation $T$  defined by 
\begin{equation}
T\eta = \psi(\eta).
\label{3.2}\end{equation}
The vector function $\psi_{n}$ is again defined by \eqref{1.14}. 
Our next lemma follows from  \eqref{2.7}.
\begin{La}\label{Lemma 3.2}The vector function $\psi_n(\eta)$  is continuous,  satisfies the relation
\begin{equation}\psi_n(\eta+q)=\psi_n(\eta)+q \quad {\mathrm{for}} \quad q\in \BZ^k,\label{3.2+}\end{equation}
 and has the form
\begin{equation}\psi_n(\eta)=\eta+a(\eta)+a[\psi(\eta)]+...+ a[\psi_{n-1}(\eta)],\label{3.3}\end{equation}
where $a(\eta)=\phi(1,\eta)-\eta$ .
\end{La}
Formula \eqref{3.3} is well known   for the case $k=1$ (see \cite[p. 104]{Arn}).
\begin{Pn}\label{Pn 3.3}   If   the  operator $T^{m}$ $($for some $m\in\BN)$  has a fixed point,
then the limit vector 
\begin{equation}\rho=\lim_{n{\to}\infty}{\psi_n(\eta)}\big/{n}\label{3.5} \end{equation}
 exists at this point and all the elements of $\rho$ are rational.
\end{Pn}
\begin{proof} 
At the fixed point $\eta$, we have
\begin{equation}T^m\eta=\psi_m(\eta)=\eta +q \quad {\mathrm{for\,\, some}} \quad q\in \BZ^k,\label{3.6} \end{equation}
because the points $\eta$ and $\eta+q$ are regarded
as identical for the torus (see Remark \ref{Remark 3.1}).
Taking into account \eqref{3.2+} and \eqref{3.6}, we obtain
\begin{align}\nn &
\psi_{2m}(\eta)=\psi_m(\eta+q)=\psi_m(\eta)+q=\eta+2q, 
\\ \nn &
 \psi_3m(\eta+q)=\psi_m(\psi_{2m}(\eta))=\eta+3q, \quad \ldots,
\end{align}
and it follows by induction that
\begin{equation}\psi_{mn}(\eta)=\eta+nq \quad (n=0,1,2,...).\label{3.7}\end{equation}
Every  $K\in \BN$ admits representation $K=m\ell+s,$ where the integer $s$ satisfies
$0 \leq s<m.$ Thus, \eqref{3.7} implies that
\begin{equation} \psi_{K}(\eta)=\psi_{s}(\eta)+q\ell.\label{3.8}\end{equation}
Using \eqref{3.8}  we derive
\begin{equation}\rho=\lim_{K\to \infty}{\psi_{K}(\eta)}\big/{K}={q}/{m},\label{3.9}\end{equation}
which proves the lemma.
\end{proof}

{\it Further we consider the general case, where fixed points \eqref{3.6} are not necessary.} Let us introduce the following notation
\begin{equation} a_{n}(\eta)= a(\eta)+a[\psi(\eta)]+...+ a[\psi_{n-1}(\eta)].\label{3.10}\end{equation}
It follows from \eqref{3.3} and \eqref{3.10} that $\psi_n(\eta)-a_n(\eta)=\eta$ so that 
\begin{equation}
 \lim_{n \to \infty}\big(\psi_n(\eta)-a_n(\eta)\big)\big/n]=0.
\label{3.13-}\end{equation} 
Let us introduce the notion of a limit set for a sequence of vectors $B_{n}$ $(n\in \BN)$ or a family of vectors $B_t$ $\big(t\in (0,\infty)\big)$.
\begin{Dn}\label{Definition 3.4} The vector $\gamma$ belongs to the  limit set  for a sequence of vectors $B_{n}$  $($for a family of vectors $B_t)$
if there is a sequence $\{n_p\}$  $(\{t_p\})$, where  $n_p\in \BN$ $\big(t_p\in (0,\infty)\big)$such
that
\begin{equation}\gamma =\lim_{n_p \to \infty}B_{n_p} \quad \big(\gamma =\lim_{t_p \to \infty}B_{t_p}\big).\label{3.12}\end{equation}\end{Dn}
Since the vector function $a(\eta)$ is continuous (together with $\psi(\eta)$) and \eqref{2.8} holds, $a(\eta)$ is bounded:
\begin{equation}\sup_{\eta\in E_k}\|{a(\eta)}\| \leq M<\infty . \label{3.13}\end{equation}
 Taking into account \eqref{3.10} and \eqref{3.13}, we see that $a_{n}(\eta)/n$ is also bounded:
\begin{equation}
\sup_{\eta\in E_k}\|{a_n(\eta)}/n\| \leq M ,
\label{3.14}\end{equation}
where $M$ is the same as in \eqref{3.13}.
We denote the limit set of  $a_{n}(\eta)/n$ by $Q(\eta).$
Relation \eqref{3.14} implies the following proposition.
\begin{Pn} \label{Proposition 3.5}The limit sets  $Q(\eta)$ of the vector functions  $a_{n}(\eta)/n$  are nonempty and uniformly bounded:
\begin{equation}
\sup_{\gamma(\eta)\in Q(\eta)}\|\gamma(\eta)\| \leq M.
\label{3.15}\end{equation}\end{Pn}
In view of \eqref{3.13-} we obtain the corollary below.
\begin{Cy}\label{Corollary 3.6} The limit set of $\psi_{n}(\eta)/n$ coincides with he limit set  $Q(\eta)$  of $a_{n}(\eta)/n$ .\end{Cy}
The set $Q(\eta)$ is the set of rotation vectors.
\begin{Tm}\label{Theorem 3.7} The limit set of $\phi(t,\eta)/t$, for $t{\to}+\infty$, coincides with $Q(\eta)$.\end{Tm}
\begin{proof} For $t\in (0,\infty)$ we choose $n_t \in \BN$ so that $ n_t \leq t< n_t+1$. Hence $t$ admits representation
\begin{equation} t=n_t+ s_t \quad (0\leq s_t<1).\label{3.16}\end{equation}
 In view of \eqref{2.5} and \eqref{3.16},  we have
\begin{equation}\phi(t,\eta)=\phi(s_t,\psi_{n_t}(\eta)).\label{3.17}\end{equation}
We represent $\psi_{n_t}(\eta)$ in the form
\begin{equation}\psi_{n_t}(\eta)=U_{n_t}(\eta)+V_{n_t}(\eta),\label{3.18}\end{equation}
where all the elements of the vector $U_{n_t}(\eta)$ are integer and all the elements of the vector $V_{n_t}(\eta)$ are nonnegative and less then 1.
Taking into account relations \eqref{2.4}, \eqref{3.17} and \eqref{3.18}, we obtain 
\begin{equation}\phi(t,\eta)=U_{n_t}(\eta)+\phi(s_t,V_{n_t}(\eta)).\label{3.19}\end{equation} 
Clearly, the norms of $V_{n_t}(\eta)$ and of  $\phi(s_t,V_{n_t}(\eta))$ for our continuous
vector function $\phi$ are bounded.
Thus, the theorem follows from Corollary \ref{Corollary 3.6} and the equalities \eqref{3.18} and \eqref{3.19}.
\end{proof}
\begin{Cy}\label{Corollary 3.8}The vector function $\phi(t,\eta)$ may be represented in the  form
\begin{equation}\phi(t,\eta)=\psi_{n_t}(\eta)+b(t,\eta),\label{3.20}\end{equation}
where $b(t,\eta)$ is a bounded vector function depending on $t\geq 0$ and $\eta \in E_k$.\end{Cy}
According to \eqref{3.20}, $\psi(\eta)$ defines  $\phi(t,\eta)$ up to a bounded term.
Using  \eqref{3.3} and \eqref{3.10},  we  rewrite  \eqref{3.20} in the form
\begin{equation}\phi(t,\eta)=\eta +a_{n_t}(\eta)+b(t,\eta),\label{3.21}\end{equation}
Taking into account \eqref{3.14} and \eqref{3.21}, we obtain the next theorem.
\begin{Tm}\label{Theorem 3.9} The vector  function 
$$F(t,\eta)=[\phi(t,\eta)-\eta]/t \quad (t>1)$$ 
is bounded and continuous. The limit set of $F(t,\eta)$, for $t{\to}+\infty$, coincides with $Q(\eta)$.\end{Tm}
Theorem \ref{Theorem 3.9} implies
 (see \cite[Ch.XYI, Theorem 1.1]{CL} or \cite[Ch.X, Assertion~7.1]{Lef}) the following proposition.
\begin{Pn}\label{Proposition 3.10} The limit set $Q(\eta)$ is closed, connected and nonempty.\end{Pn}
Note that the fact that $Q(\eta)$ is closed and nonempty easily follows from the considerations above
Theorem \ref{Theorem 3.9}.

Let $q$ be a vector with integer elements. Then, taking into account \eqref{2.8} \eqref{3.2+} and \eqref{3.10} we  have
\begin{equation}a_n(\eta+q)=a_n(\eta),\quad n=1,2,...\label{3.22}\end{equation}
It follows from \eqref{3.22} that
\begin{equation}Q(\eta+q)=Q(\eta).\label{3.23}\end{equation}
\section{ Jacobian matrix, Zorn's lemma}
Suppose that a $k$-vector function f of $k$ variables is differentiable.
 Then, the corresponding Jacobian matrix has the form
\begin{equation}J_{x}(f)=\left(
                               \begin{array}{ccc}
                                 \frac{\partial{f_1}}{\partial{x_1}}&... & \frac{\partial{f_1}}{\partial{x_k}} \\
                                 ... & ... &... \\
                                 \frac{\partial{f_n}}{\partial{x_1}} & ... & \frac{\partial{f_k}}{\partial{x_k}} \\
                               \end{array}
                             \right)\label{4.1}\end{equation}
 Differentiable $k$-vector functions $f(x)$ and $g(f(x))$
 satisfy the chain rule:
\begin{equation}J_{x}(g(f(x))=J_{f}(g)J_{x}(f).\label{4.2}\end{equation}
\emph{In this section, we assume that all the first  order derivatives of k-vector function $\psi(\eta)$, $(\eta{\in}E_{k})$ exist and are continuous. Then,  all the first  order derivatives of $k$-vector functions $\psi_n(\eta)$ $(n\in \BN)$ exist and are continuous as well.}

Here, we do not require that the vector function $\psi(\eta)$ necessarily satisfies condition  \eqref{2.4}, that is, the condition $ \psi(\eta+q)=\psi(\eta)+q$.
It follows from \eqref{4.2} that  the Jacobian for $\psi_n(\eta)$ has  the form:
\begin{equation}J_{\eta}(\psi_{n}(\eta))=J_{\psi_{n-1}}(\psi_{n})J_{\psi_{n-2}}(\psi_{n-1})...J_{\eta}(\psi),\label{4.3}\end{equation}
where $\psi_{0}(\eta)=\eta$.

\begin{Cy}\label{  Corollary 4.1} If  $\det[J_{\eta}(\psi)]$ has no roots, then  $\det[J_{\eta}(\psi_{n})]$ has no roots too.\end{Cy}

\begin{Cy}\label{Corollary 4.2}If  $\det[J_{\eta}(\psi)]$ has  no roots,  then the vector  functions $\psi_n(\eta)$ $(n\in \BN)$  have no singular points.\end{Cy}
Here,   $\eta$ is called a  singular point for the function  $\psi(\eta)$  if
$\det[J_{\eta}(\psi)]=0$.
\begin{Cy}\label{  Corollary 4.3} If  $\det[J_{\eta}(\psi)]$ has a root in the point $\eta=X$, then  $\det[J_{\eta}(\psi_{n})]$ has  a root in the point $\eta=X$  too.\end{Cy}

\begin{Cy}\label{Corollary 4.4}If  $\det[J_{\eta}(\psi)]$ has a root in the point $\eta=X$,
  then the point $X$ is a singular point of the  vector  functions $\psi_n(\eta)$ $(n\in \BN)$.\end{Cy}
Next, we will need the well-known Zorn's lemma. First, we  introduce the corresponding notions of the maximal and minimal elements.
\begin{Dn}\label{Definition 4.6}A vector $X$ in  a semi-ordered set $S$ is called maximal $($minimal$)$, if there is no other element of $S$ greater $($lesser$)$
than $X$. \end{Dn}
\begin{La}\label{Lemma 4.5}\cite{Zor} A semi-ordered set $S$ containing upper  bounds for every totally ordered subset contains at least one maximal element.\end{La}
Zorn's lemma above is easily reformulated in terms of the minimal elements. The points $\eta$ of a semi-ordered set $S$, where $\psi(\eta)$ takes its maximal or minimal
(for the set $\{\psi(\eta): \, \eta \in S\}$) values, are called extremal points.
\begin{Tm}\label{Theorem 4.7}
 Let all the first  order derivatives of the $k$-vector function $\psi(\eta)$ $(\eta{\in}S\subset E_{k})$ exist and be continuous. If the subset $S$ of $E_k$ is closed and bounded
and the inequality  
\begin{equation}\det[J_{\eta}(\psi)]{\ne}0 \quad (\eta \in S)
\label{4.4}\end{equation}
is valid, then the following statements hold.

1.  The set of extremal points of the $k$-vector function $\psi(\eta)$, where $\eta$ belongs to $S\in E_k$, is not empty.

2. All the extremal points of the $k$-vector function $\psi(\eta)$ $(\eta{\in}S)$ belong to the boundary of the subset $S$.
\end{Tm}
\begin{proof}
The assertion 1 of the theorem follows from the Zorn's lemma.
Let us prove the assertion 2.
According to condition  \eqref{4.4}, the vector function  $\psi(\eta)$ is locally bijective.  Hence, for any interior point $X\in S$, the vector
$\psi(X)$  has in its neighbourhood a point Y such that $\psi(X)$ is less than $Y$.
The assertion~2 is proved,
\end{proof}
\begin{Rk}\label{Remark 4.8} Theorem 4.7 is a generalisation to the vector case of the well-known theorem for the scalar case:\\
If a function $\psi$ is continuous on a closed and bounded set $S{\in}E_k$ and takes its greatest (smallest) value at  $X\in S$,
then either $X$ belongs to the interior of  $S$ and $\psi$ takes a local    greatest (smallest) value at $X$ or $X$ lies on the boundary set of $S$.\end{Rk}
Further we assume that subset $S$ of the set $E_k$ is closed and bounded and the function $\psi(\eta)$ is continuous. Then the set of values of $\psi(\eta),\,\eta{\in}S$
is also closed and bounded. Therefore, according to  Zorn's lemma there exist maximal and minimal values of $\psi(\eta)$ $(\eta{\in}S)$.

\emph{Using the  maximal and minimal values of $\psi(\eta)$ we  construct the smallest vector $X{\in}E_k$ and the greatest vector $Y{\in}E_k$ such that}
\begin{equation}X\geq \psi(\eta)\geq Y \quad (\eta \in S). \label {4.5}\end{equation}
For this purpose, we consider extremal points of $\psi(\eta)$. We denote the set of all values of $\psi(\eta)$ ($\eta \in S$) by $\Psi$ and 
the set of all maximal (minimal) values of $\Psi$    by $M$ ($L$). The sets of $p$th elements $\psi_p(\eta)$  of $\psi(\eta){\in}M$ and $\psi(\eta){\in}L$ are
denoted by $M_p$ and $L_p$, respectively.
\begin{Tm}\label{Theorem 4.9} Let the set $S\subset E_k$ be closed and bounded and let the $k$-vector function
$\psi(\eta)$ be continuous on $S$. Let the vectors $X$ and $Y$ be given by the equalities
\begin{equation} X=[m_1,m_2...,m_k] ,\quad Y=[\ell_1,\ell_2...,\ell_k],\label{4.6}\end{equation}
where 
\begin{equation}
m_p=\sup_{z_p\in M_p}z_p,\quad \ell_p=\inf_{z_p\in L_p}z_p \quad (1\leq p \leq k).
\label{4.7}\end{equation}
Then, \eqref{4.5} holds for $X$ and $Y$ introduced above. Moreover, this $X$ is the smallest and $Y$ is the  greatest of the vectors satisfying \eqref{4.5}.

\end{Tm}
\begin{proof} Let us  prove \eqref{4.5} by negation. If \eqref{4.5} does not hold, then either a set $\Psi_{\max}$ of $\psi(\eta)$ such that 
$\psi_p(\eta) \geq m_p+\varepsilon$ or a  set  $\Psi_{\min}$ of $\psi(\eta)$ such that 
$\psi_p(\eta) \leq \ell_p-\varepsilon$ is nonempty, bounded and closed in $E_k$ for some $p$ $(1\leq p\leq k)$ and $\varepsilon>0$.
According to Zorn's  lemma, there is a maximal (minimal) element of $\Psi_{\max}$ ($\Psi_{\min}$). If  $\Psi_{\max}$ is nonempty,
its maximal elements are also maximal in $\Psi$ since there is not a vector in $\Psi - \Psi_{\max}$, which is greater than one of the vectors
from $\Psi_{\max}$. At the same time, the maximal vectors from $\Psi_{\max}$ do not belong $M$ and we arrive at a contradiction.
In a similar way, we arrive at a contradiction if $\Psi_{\min}$ is nonempty. Thus, \eqref{4.5} is proved.

The assertion that $X$ is the smallest and $Y$ is the  greatest of the vectors satisfying \eqref{4.5} easily follows from the construction
of $X$ and $Y$  via \eqref{4.6} and \eqref{4.7}.
\end{proof}

 Relations  \eqref{4.6} and \eqref{4.7}  imply the following corollary.
 \begin{Cy}\label{Corollary 4.11}If the set $M$ contains only one maximal (minimal) vector, then this vector is $X$ ($Y$). \end{Cy}
 \section{Examples}
\begin{Ee}\label{Example 5.1}  Starting with the simplest case, where  relations \eqref{2.1} and \eqref{2.2} hold, let
\begin{equation}\frac{dx}{dt}=G,\quad x(0)=\eta \qquad (t\in \BR; \,\, \eta,\,G,\,x(t)\in E_k). \label{5.1}\end{equation}
\end{Ee}
The solution $\phi(t,\eta)$ of equation \eqref{5.1} and $\psi(\eta)$ have the form
\begin{equation} \phi(t,\eta)=Gt+\eta, \quad \psi(\eta)=G+\eta .\label{5.2}\end{equation}
It follows that
\begin{equation} \psi_n(\eta)=nG+\eta,\quad Q(\eta)=\{G\} \quad ({\mathrm{i.e.,}}\, \rho_1=\rho_2=G).
\label{5.3}\end{equation}
The next assertion is easily checked directly.
\begin{Pn}\label{Proposition 5.2} In the case \eqref{5.1}, relations \eqref{2.3}--\eqref{2.5} are fulfilled.
The vector functions $\psi(\eta)$ and $\phi(t,\eta)$ $($for each $t)$ are injective. \end{Pn}
\begin{Rk}\label{Remark 5.3} 
We note that the injective property of the vector function $\psi(\eta)$ is an analog of the monotonic property of $\psi(\eta)$ in the scalar
case.
Injective property of the vector function $\phi(t,\eta)$ is an  analog of  Condition II on $\phi(t,\eta)$ in the scalar case.
Thus, the injective properties  of $\psi(\eta)$ and $\phi(t,\eta)$ are of interest.
\end{Rk}
Next, consider a more non-trivial example.
\begin{Ee}\label{Example 5.4} Let $k=2$, $\eta=[\eta_1, \, \eta_2]$ and $\psi(\eta)$ have the form
\begin{equation}\psi(\eta)=[\eta_1+r\sin(2\pi\eta_2), \, \eta_2-r\sin(2\pi\eta_1)] \quad (r\in \BR).\label{5.4}\end{equation}
Let the 2-vector function $\Phi(t,\eta)$ $(0\leq t \leq 1)$ have the form
\begin{align} &
\Phi(t,\eta)=[\Phi_1(t,\eta), \, \Phi_2(t,\eta)],\quad
\Phi_1(t,\eta)=\eta_1+r\sin(\pi{t}/2)\sin(2\pi\eta_2),\,\label{5.5}
\\ &
\Phi_2(t,\eta)=\eta_2-r\sin(\pi{t}/2)\sin(2\pi\eta_1).\label{5.6}\end{align}
\end{Ee}
It is easy to see that conditions \eqref{2.10} and \eqref{2.11} are fulfilled and we may apply Theorem \ref{Theorem 2.3}.
(Clearly, the second equality in \eqref{2.4} holds for $\psi$ as well.)
\begin{Pn}\label{Proposition 5.5}
Given \eqref{5.4}--\eqref{5.6} we recover solution $\phi(t,\eta)$ of  a generalised system using  \eqref{2.12}.
\end{Pn}
Let us consider $\psi(\eta)$ in greater detail. The following inequality is well known:
\begin{equation} |\sin{x}-\sin{y}|{\leq}|x-y|.\label{5.7}\end{equation}
Using \eqref{5.4} and \eqref{5.7}, we derive
\begin{equation}\psi(\eta)-\psi(p)=\eta-p+rb,\quad \norm{b}<2\pi\norm{\eta-p},\label{5.8}\end{equation}
where $b=[\sin(2\pi\eta_2)-\sin(2\pi{p_2}), \, \sin(2\pi{p_1})-\sin(2\pi\eta_1)]$.
Taking into account \eqref{5.8}, we obtain
\begin{equation}\norm{\psi(\eta)-\psi(p)} \ne 0 \quad {\mathrm{for}} \quad \eta \ne p \quad {\mathrm{and}} \quad |r|<1/(2\pi).\label{5.9}\end{equation}
Using relation \eqref{5.9} we have:
\begin{Pn}\label{Proposition 5.6} The vector function $\psi(\eta)$ given by \eqref{5.4}, where $\break |r|<1/(2\pi)$, is an injective mapping of the  space  $E_2$ into the space $E_2$.
\end{Pn}
\begin{Rk}\label{Remark 5.7} Proposition \ref{Proposition 5.6}
  may be generalised using representation \eqref{2.7} of $\psi$ in the $k$-dimensional generalised system. Namely, if
the inequality
\begin{equation}\norm{a(\eta)-a(p)}<\norm{\eta-p}.\label{5.10}\end{equation}
is fulfilled for all the $k$-vectors $\eta$ and $p$, then the vector function $\psi(\eta)$ is an injective mapping of the  space  $E_k$ into the space $E_k$.
It is easy to see that the vector function $\psi(\eta)$ considered in Proposition \ref{Proposition 5.6} satisfies the condition \eqref{5.10}. \end{Rk}

Now, let us calculate the Jacobian matrix $\textbf{J}$ for this $\psi(\eta)$:
\begin{equation}\textbf{J}=\left(
                  \begin{array}{cc}
                    1 & 2r\pi\cos(2\pi\eta_2) \\
                    -2r\pi\cos(2\pi\eta_1) & 1 \\
                  \end{array}
                \right).\label{5.11}\end{equation}
Hence, we have
\begin{equation} \det\textbf{J}=1+(2{\pi}r)^2\cos(2\pi\eta_2)\cos(2\pi\eta_1)>0 \quad {\mathrm{for}} \quad |r|<1/(2\pi).\label{5.12}\end{equation}
Formula \eqref{5.12} yields the next proposition.
\begin{Pn}\label{Proposition 5.8}
The vector function $\psi(\eta)$ given by \eqref{5.4}, where $\break |r|<1/(2\pi)$, has no critical points.\end{Pn}
Recall  that a \emph{critical point} of the vector function is a point where the rank of the Jacobian matrix is not maximal.

Similar to Proposition \ref{Proposition 5.6} we derive a more general proposition.
\begin{Pn}\label{Proposition 5.9} The vector function $\Phi(t,\eta)$ given by \eqref{5.5} and \eqref{5.6},
where $|r|<1/(2\pi)$,  is an injective mapping of the  space  $E_k$ into the space $E_k$ for each $t\in [0,1]$.
\end{Pn}
For the Jacobian matrix ${\bf J}=\{{\bf J}_{ik}\}$  $\,(1\leq i \leq 2, \,\,  1\leq k \leq 3)$ of  $\Phi(t,\eta)$ we easily obtain:
$ {\bf J}_{11}={\bf J}_{22}=1$,
$$ {\bf J}_{12}=2r\pi\sin(\pi{t}/2)\cos(2\pi\eta_2), \quad {\bf J}_{21}=-2r\pi\sin(\pi{t}/2)\cos(2\pi\eta_1).$$
Thus, the rows of  ${\bf J}$ are linearly independent for $|r|<1/(2\pi)$ and    the Jacobian matrix $\textbf{J}$ has the maximal rank $2$. 
\begin{Pn}\label{Proposition 5.10}
The vector function $\Phi(t,\eta)$  given by \eqref{5.5} and \eqref{5.6},
where $|r|<1/(2\pi)$, has no critical points.\end{Pn}
\begin{Ee}\label {Example 5.11} Let us consider the vector function $\psi(\eta)$ of the form \eqref{5.4}, where $r=1/{2\pi}$.\end{Ee}

It follows from \eqref{5.12} that
\begin{equation} \det{J(\eta)}=0 \quad (\eta=[\eta_1,\eta_2])\label{5.13}\end{equation}
for the cases
\begin{align}& \eta_1=p,\quad \eta_2=1/2+s \quad (p,s \in \BZ);\label{5.15}
\\ &
\eta_2=p, \quad \eta_1=1/2+s \quad (p,s \in \BZ).\label{5.16}\end{align}
Hence, we have:
\begin{Cy}\label{Corollary 5.12}The points  \eqref{5.15} and \eqref{5.16} are singular points of the  vector functions $\psi_n(\eta).$ \end{Cy}




{\bf Acknowledgement.} The author is very grateful to Alexander Sakhnovich for the careful reading of the work and many useful remarks.


\begin{thebibliography}{11}
\bibitem{ARS}
Aliste-Prieto, J.,  Rand, B., Sadun, L.:
Rotation numbers and rotation classes on one-dimensional tiling spaces.  
Ann. Henri Poincar\'e {\bf 22}, 2161--2193 (2021) 
\bibitem{Arn}
 Arnold, V.I.: { Geometrical methods in the theory of ordinary differential equations.}  Springer, New York (1988)
\bibitem{Bo}
Bohl, P.: { \"{U}ber die Hinsichtlich der Unabh\"{a}ngigen und Abh\"{a}ngigen Variabeln Periodische Differentialgleichung Erster Ordnung. } Acta Math. {\bf 40}, 321--336 (1916)
\bibitem{CL}
Coddington, E.A., Levinson, N.:  {Theory of ordinary differential equations.} McGraw-Hill Book Co., New York  (1955)
\bibitem{Fey}
Feynman, R.: {QED: The strange theory of light and matter.} Prinston University Press, Prinston (1985)
\bibitem{Fra1}
Franks, J.: Realizing rotation vectors for torus homeomorphisms. Trans. Amer. Math. Soc. {\bf 311}, 107--115 (1989)
\bibitem{Lef}
Lefschetz, S.: Differential  equations: geometric theory. Interscience Publishers, New York (1957)
\bibitem{Lel}
Lellouch, G.: On rotation sets for surface homeomorphisms in genus $\geq 2$. (French)  Mem. Soc. Math. Fr. (N.S.) No. 178 (2023)
\bibitem{Lli0}
Llibre, J., MacKay, R. S.:
Rotation vectors and entropy for homeomorphisms of the torus isotopic to the identity. 
Ergodic Theory Dynam. Systems {\bf 11},  115--128 (1991) 
\bibitem{Lli}
Llibre, J., Valls, C.: Darboux theory of integrability in $\BT^n$. Eur. J. Math. {\bf 10},  Paper No. 26 (2024) 
\bibitem{Poi}
Poincar\'e, H.:  Sur les courbes definies par une equation differentielle. Oeuvres, vol. 1, Paris (1892)
\bibitem{Sak}
Sakhnovich, L.A.:  {The deviation factor and divergences in quantum electrodynamics, concrete examples.}  Phys. Lett. A  {\bf 383}, 2846--2853 (2019)
\bibitem{Vul}
Vulikh, B.Z.:  {Introduction to the theory  of partially ordered spaces.} Wolters-Noordhoff Scientific Publications, Gr\"oningen (1967)
\bibitem{Zor}
Zorn, M.:  {A remark on method in transfinite  algebra.} Bull. Amer. Math. Soc. {\bf 41}, 667--670
(1935)
\end{thebibliography}
\end{document}